\documentclass{amsart}
\theoremstyle{plain}
\newtheorem{theorem}{Theorem}

\newtheorem{lemma}{Lemma}

\theoremstyle{definition}

\usepackage{graphicx}

\theoremstyle{remark}

\newcommand \gluplus {\mathop{+_{\textup{glu}}}}

\DeclareMathOperator \Sub {Sub}

\newcommand \nplu {\mathbb N^+}
\newcommand \NS {NS}
\newcommand \nsub [1] {\textup{\NS}(#1)}

\newcommand \badgood[2]{{#2}}
\newcommand\nonparallel {\mathrel{\not{\kern-1pt{\mathord\parallel}}}}
\newcommand\set [1]{\{#1\}}
\usepackage{enumitem}
\usepackage{amsmath}

\numberwithin{equation}{section}

\allowdisplaybreaks[4]

\begin{document}
\title[Subuniverses of finite semilattices]{Analyzing Subuniverse Counts in Finite Semilattices: Unveiling the Rankings and Descriptions} 

\author{Delbrin Ahmed}
\address{University of Duhok\\University campus \\Zakho street \\
  Duhok\\Kurdistan region-Iraq}

\email{delbrin.ahmed@uod.ac}
\urladdr{
https://www.uod.ac/ac/c/cbe/departments/mathematics/academic-members/delbrin-ahmed/}

\author{Muwafaq Salih}

\address{University of Duhok\\University campus \\Zakho street \\
  Duhok\\Kurdistan region-Iraq}
  
\email{muwafaq.salih@uod.ac}
\urladdr{https://www.uod.ac/ac/c/cbe/departments/mathematics/academic-members/muwafaq-mahdi-salih/}

\author{Dilbak Haje }
\address{University of Duhok\\University campus \\Zakho street \\
  Duhok\\Kurdistan region-Iraq}
  
\email{dilpakhaje@uod.ac }

\begin{abstract}
Let $(L,\vee)$ be a finite n-element semilattice where $n\geq 5$. We prove that the fourth largest number of subuniverses of an $n$-element semilattice is $25\cdot 2^{n-5}$, the fifth largest number is  $ 24.5\cdot 2^{n-5}$, and the sixth one is   $ 24\cdot 2^{n-5}$.  Also, we describe the $n$-element semilattices with exactly  $25\cdot 2^{n-5}$,  $ 24.5\cdot 2^{n-5}$ or  $ 24\cdot 2^{n-5}$ subuniverses.

\end{abstract}



\maketitle

\section{Introduction and our result}
Throughout this paper, all semilattices are assumed to be finite. For  a semilattice  $(L,\vee)$, let $\Sub(L,\vee)$  denote its \textit{subuniverse-lattice}.  A subuniverse in this context refers to either a subsemilattice or the emptyset. Our notation and terminology are standard, see, for example, Chajda et al \cite{chajda2007semilattice}. However, we provide a brief review of some notations and introduce additional auxiliary concepts for clarity.

\begin{figure}[htp] 
	\centerline
	{\includegraphics[scale=0.6]{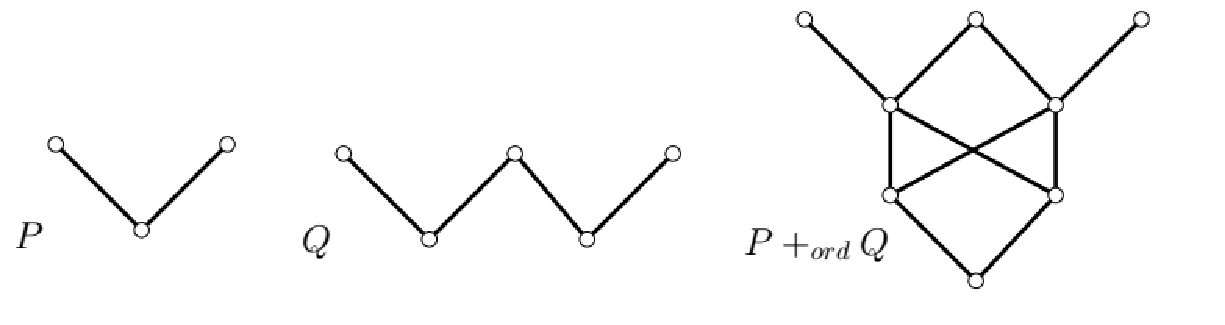}}
	\caption{{The ordinal sum $P +_{ord} Q$ of $P$ and $Q$ }
		\label{ord}}
\end{figure}

\begin{figure}[htp] 
	\centerline
	{\includegraphics[scale=0.7]{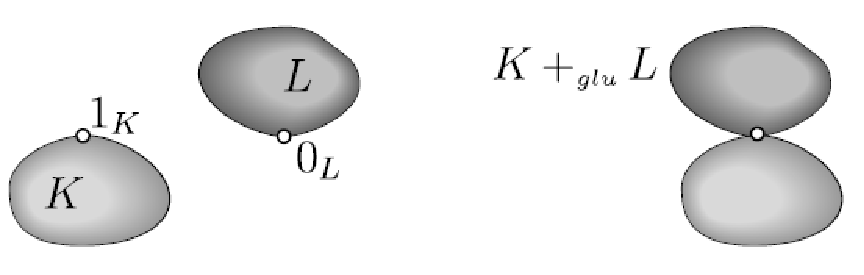}}
	\caption{{The glued sum $K \gluplus L$ of $K$ and $L$ }
		\label{glu}}
\end{figure}

On semilattice $(L,\leq)$, we have a natural partial ordering of a semilattice is defined by
$$x \leq y \iff x\lor y =y.$$
Conversely, if $(L,\leq)$ is partial order in which any two elements $x,y$ have a least upper bound $x\lor y$, then $(L,\lor)$ is a semilattice. 
For any $x,y$ in a join-semilattice, $x \land y$  is defined by their infimum provided it  exists; if this infimum does not exist, then $ x\land y$  is undefined. Let $P$ and $Q$ be posets with disjoint underlying sets. Then the \textit{ordinal sum} $P +_{ord} Q$ is the poset on $P \cup Q$ with $s \leq t$ if either $s,t \in P$ and $s \leq t $; or $s, t \in Q$ and $ s\leq t$; or $s \in P$ and $t \in Q$, 
To draw the Hasse diagram of $P +_{ord} Q$, we place the Hasse diagram of $Q$ above that of $P$ and then connect any minimal element of $Q$ with any maximal element of $P$; see Figure \ref{ord}. If $K$ with 1 and $L$ with 0 are finite posets, then their glued sum $K \gluplus L$ is the ordinal sum of the posets $K \setminus \{1_K$\}, the singleton poset, and $L \setminus \{0_L\}$, in this order; see Figure \ref{glu}. 

We inspired by similar or analogous results concerning lattices and semilattices , see Ahmed and Horv\'ath \cite{ahmed_horvath_2019, ahmed2021first}, Cz\'edli  \cite{oaumb2018003} \cite{czedli2019eighty} \cite{czedli2019finite}  \cite{czedli2019lattices}  \cite{czedli2019one} and Cz\'edli and Horv\'ath \cite{czedli2018note}.

For a more in-depth exploration of lattice theory and semilattices,  we direct the reader to  the references provided in \cite{gratzer2011lattice} \cite{rival1979lattices} and  \cite{chajda2007semilattice} for detailed insights into these subjects.

For a natural number $n\in\nplu:=\set{1,2,3,\dots}$, let 
\begin{equation*}
\nsub n:=\set{|\Sub(L,\vee)|: (L,\vee)\textup{ is a semilattice of size }|(L,\vee)|=n}.
\end{equation*}
\newpage

\begin{theorem}\label{thmmain}
	If $5\leq n\in \nplu $, then the following assertions hold.
	\begin{enumerate}[label=(\roman*)]
	
		\item\label{thmmaind}  The fourth largest number in $ \nsub n $ is  $  25\cdot 2^{n-5} $. Furthermore, an $n$-element semillatice $(L,\vee)$ has exactly $ 25\cdot 2^{n-5 }$ subuniverses if and only if $(L,\vee)\cong  H_5 \gluplus C_1 $ or $(L,\vee)\cong C_0 +_{ord} H_5 \gluplus C_1 $, where \badgood{$C_0$ , $C_1$ and  $C_2$ }{$C_0$ and $C_1$} are finite chains.
		
		\item\label{thmmaine}  The fifth largest number in $\Sub(L,\vee)$ is $ 24.5\cdot 2^{n-5}.$ Furthermore, an $n$-element semilattice $(L,\vee)$ has exactly $  24.5\cdot 2^{n-5} $ subuniverses if and only if $(L,\vee)\cong H_3\gluplus B_4 \gluplus C_1$ or
		$(L,\vee)\cong C_0 +_{ord} H_3\gluplus B_4 \gluplus C_1$, where \badgood{$C_0$ , $C_1$ and  $C_2$ }{$C_0$ and $C_1$ } are finite chains.

		\item\label{thmmainf}  The sixth largest number in $\Sub(L,\vee)$ is $ 24\cdot 2^{n-5}.$ Furthermore, an $n$-element semilattice $(L,\vee)$ has exactly $  24\cdot 2^{n-5} $ subuniverses if and only if $(L,\vee)\cong K_3 \gluplus C_1$ or
		$(L,\vee)\cong C_0 +_{ord} K_3 \gluplus C_1$, where \badgood{$C_0$ , $C_1$ and  $C_2$ }{$C_0$ and $C_1$ } are finite chains.

	\end{enumerate} 
\end{theorem}

\section{Preparatory lemmas}\label{sectionprepare}

An element $u$ of a semilattice $L$ is called a narrow element, or \it a narrows \rm for short,
if  $u \neq 1_L$ and
$L ={\uparrow} u \cup {\downarrow}u$. That is if $u \neq 1_L$ and   $x \| u$ holds for no $x \in  L$.

The concept of a binary partial algebra is widely recognized. However, for a refresher, readers can refer to \cite{czedli2019eighty}.  Let $\mathcal{A}$ be a finite $n$-element binary partial algebra. A \it subuniverse \rm of $\mathcal{A}$ is a subset $X$ of $\mathcal{A}$ such that $X$ is closed with respect to all partial operations. The set of subuniverses of $\mathcal{A}$ will be denoted by $\Sub(\mathcal{A}).$  The \it relative number of subuniverses of  $\mathcal{A}$ \rm denoted by  $ \sigma_k(\mathcal{A}) $ is defined as follows:
$$ \sigma_k(\mathcal{A}) = | \Sub(\mathcal{A})|\cdot 2^{k-n}.$$

\noindent In this paper we will use $k=5$. An original definition of $\sigma_k$ is given in the paper of Cz\'edli  \cite{czedli2019eighty}, where he used $k=8$. 

\begin{lemma}\label{lemma2.1-7}\cite{ahmed2021first}
	If $(K,\vee)$ is a subsemilattice and $H$ is a subset of a finite semilattice $(L,\vee)$,
	then the following three assertions hold.
	\begin{enumerate}[label=(\roman*)]
		\item\label{lemm21a}With the notation $t:= |{ H \cap S : S \in Sub(L,\vee) }|$, we have that
		$$ \sigma_k(L,\vee) \leq t \cdot 2^{k-|{H|} }. $$
		\item\label{lemm21b}$\sigma_k(L,\vee) \leq \sigma_k(K,\vee).$
		\item\label{lemm21c} Assume, in addition, that $(K,\vee)$ has no narrows. Then $ \sigma_k(L,\vee) =  \sigma_k(K,\vee)$ if and only if $(L,\vee)$ is (isomorphic to) $C_0 +_{ord}(K,\vee) \gluplus C_1,$ where $C_1$ is a chain, and $C_0$ is a chain or the emptyset.  
	\end{enumerate}
	
\end{lemma}


The following lemma can be proved by using a computer program,   a computer program is accessible on  G. Cz\'edli  webpage:  http://www.math.u-szeged.hu/\textasciitilde{}czedli/.

\begin{figure}[htb] 
	\centerline
	{\includegraphics[scale=1]{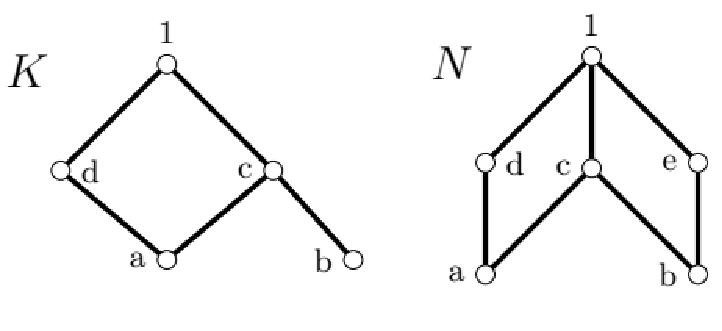}}
	\caption{{ Partial lattices $K_3$, $K$ and $N$}
		\label{figsix}}
\end{figure}%

\begin{figure}[htb] 
	\centerline
	{\includegraphics[scale=1]{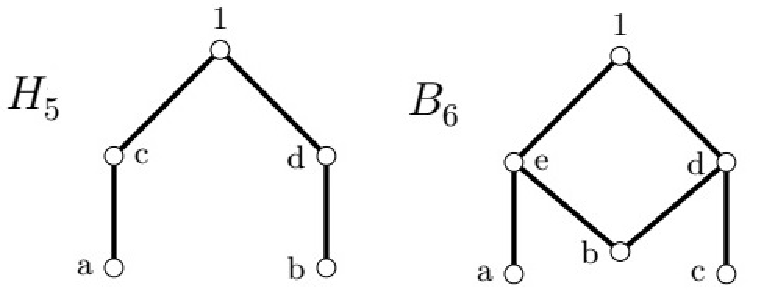}}
	\caption{{Partial lattices $ H_5$ and  $K_0$}
		\label{figseven}}
\end{figure}%

\begin{figure}[htb] 
	\centerline
	{\includegraphics[scale=1]{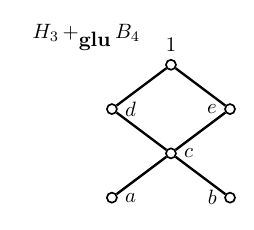}}
	\caption{{Partial lattice $ H_3\gluplus B_4$}
		\label{figeight}}
\end{figure}%

\begin{lemma}\label{lemmaLHFtfsZh}
	For the join-semilattices given in figures~\ref{figsix} to  \ref{figeight}  the following six assertions hold.
	\begin{enumerate}[label=(\roman*)]
		
		\item\label{iii}$\sigma_5(H_5)=25,$
		\item\label{ii}$ \sigma_5(H_3\gluplus B_4)=24.5,$                                                                          
		\item\label{iiii}$\sigma_5(K_3)=24,$
		\item\label{v}$  \sigma_5(K)=23,$                                        
		\item\label{ix}$  \sigma_5(N)=19.5,$                                      
		\item\label{xi}$  \sigma_5(K_0)=15.25$.
		
	\end{enumerate}
\end{lemma}

\begin{proof} 
	
	The notations  given by figure~\ref{figsix} to ~\ref{figeight} will  be used. For  later reference, note that if $ (L,\vee) $ is a chain then $|Sub(L,\vee)|= 2^{|(L,\vee)|}$.
	

	For \ref{iii}, let us compute
	\begin{align*}
	&|\set{S\in \Sub(H_5,\vee): d\not\in S}|=\badgood{9}16,\qquad\text{ (S is chain),}\cr
	&|\set{S\in \Sub(H_5,\vee): d\in S, \{a,b,c\}\cap\ S =\emptyset}|=2, \text{ and}\cr
	&|\set{S\in \Sub(H_5,\vee): d\in S, \{a,b,c\}\cap\ S \neq\emptyset}|=7.
	\end{align*}
	Hence, $|\Sub(H_5,\vee)|=\badgood{9}16+2+7=\badgood{24}{25}=25\cdot 2^{5-5}$, which means that  $\sigma_5(H_5)=25$; and this  proves case \ref{iii}.

	For \eqref{ii}, let us compute
	\begin{align*}
	&|\set{S\in \Sub(H_3\gluplus B_4,\vee): a\not\in S}|=\badgood{9}28,\cr
	&|\set{S\in \Sub(H_3\gluplus B_4,\vee): a\in S, \{b,c,d,e\}\cap\ S =\emptyset}|=2, \text{ and}\cr
	&|\set{S\in \Sub(H_3\gluplus B_4,\vee): a\in S, \{b,c,d,e\}\cap\ S \neq\emptyset}|=19.
	\end{align*}
	Hence, $|\Sub(H_3\gluplus B_4,\vee)|=\badgood{9}28+2+19=\badgood{24}{49}=24.5\cdot 2^{6-5}$, while  $\sigma_5(H_3\gluplus B_4)=24.5$ proves \eqref{ii}.

	For \ref{iiii}, let us compute
	\begin{align*}
	&|\set{S\in \Sub(K_3,\vee): a\not\in S}|=\badgood{9}7,\qquad\text{ (by Lemma 3.3.2 (i) in \cite{ahmed2021first}),}\cr
	&|\set{S\in \Sub(K_3,\vee): a\in S, \{b,c\}\cap\ S =\emptyset}|=2, \text{ and}\cr
	&|\set{S\in \Sub(K_3,\vee): a\in S, \{b,c\}\cap\ S \neq\emptyset}|=3.
	\end{align*}
	Hence, $|\Sub(K_3,\vee)|=\badgood{9}7+2+3=\badgood{24}{12}=24\cdot 2^{4-5}$, which means that   $\sigma_5(K_3)=24$; and this  proves case \ref{iiii}.

	In order to prove \ref{v}, note that $(B_4,\vee)$  has $|\Sub(B_4)|+1$.
	Now, let us compute
	\begin{align*}
	&|\set{S\in \Sub(K,\vee): b\not\in S}|=\badgood{9}14,\qquad\text{($S$ is $B_4$),}\cr
	&|\set{S\in \Sub(K,\vee): b\in S, \set{a,c,d,}\cap S=\emptyset}|=2, \text{ and}\cr
	&|\set{S\in \Sub(K,\vee): b\in S,\set{a,c,d,}\cap S\neq\emptyset}|=7,
	\end{align*}
	whereby $|\Sub(K,\vee)|=\badgood{9}14+2+7=\badgood{24}{23}=23\cdot 2^{5-5}$, which means that   $\sigma_5(K,\vee)=23$; and this  proves case  \ref{v}.  
	

	
	
	For \ref{ix}, let us compute
	\begin{align*}
	&|\set{S\in \Sub(N,\vee): d \not\in S}|=\badgood{9}23,\qquad\text{(by \ref{v}),}\cr
	&|\set{S\in \Sub(N,\vee): d\in S , \set{a,b,c,e}\cap S=\emptyset}|=2, \text{ and}\cr
	&|\set{S\in \Sub(N,\vee): d\in S,\set{a,b,c,e}\cap S\neq\emptyset}|=14.
	\end{align*}
	Hence, $|\Sub(N,\vee)|=\badgood{9}23+2+14=\badgood{24}{39}=19.5\cdot 2^{6-5}$, which means that   $\sigma_5(N,\vee)=19.5$; and this  proves case \ref{ix}. 
	
	For \ref{xi}, let us compute
	\begin{align*}
	&|\set{S\in \Sub(K_0,\vee): b \not\in S}|=\badgood{9}39,\qquad\text{(by \ref{ix}),}\cr
	&|\set{S\in \Sub(K_0,\vee): b\in S , \set{a,x,c,y,z}\cap S=\emptyset}|=2, \text{ and}\cr
	&|\set{S\in \Sub(K_0,\vee): b\in S,\set{a,x,c,y,z}\cap S\neq\emptyset}|=20.
	\end{align*}
	Hence, $|\Sub(K_0,\vee)|=\badgood{9}39+2+20=\badgood{24}{61}=15.25\cdot 2^{7-5}$, which means that   $\sigma_5(K_0,\vee)=15.25$; and this  proves case \ref{xi}.

\end{proof}

\section{The rest of the proof}

\begin{proof}[Proof of Theorem~\ref{thmmain}]

We will prove part \ref{thmmaind}.

	Assume that $ (L,\vee)$ is of the given form. Then $\sigma_5(L,\vee)=25$ is clear from Lemma \ref{lemma2.1-7}\ref{lemm21c} and Lemma \ref{lemmaLHFtfsZh}\ref{iii}. In order to prove the converse, assume that $\sigma_5(L,\vee)=25$. By Theorem 1 (iii) in \cite{ahmed2021first},
	$(L,\vee)$ has two incomparable elements, $a$ and $d$.  From Theorem 1 (iii) in \cite{ahmed2021first}, $\{a,d\}$
is not the only 2-element antichain in $L$ since otherwise $\sigma_5(L,\vee)$  would equal 28. Also, if we have two 2-element antichain in $L$ then $\sigma_5(L,\vee)=26$ by Theorem \ref{thmmain}(iii)  in \cite{ahmed2021first}. Therefore we know $L$ has more than two antichains. So we have $\{a,b\}$, $\{c,d\}$ and $\{e,f\}$, where the elements $a, b, c, d, e, f $ are distinct. Let $x= a \vee b$, $y= c \vee d$ and $z= e \vee f$. There are cases that depend on  $t:=| \{a,b,c,d,e,f,x,y,z\}|,$ which is 9,8,7 or 6. The number of cases is reduced by the symmetry roles of $ a, b, y, z$. Now the number of elements is $t=9$. Take the partial algebra $U_1=\{a,b,c,d,e,f,x,y,z\}$ with  $x= a \vee b$, $y= c \vee d$ and $z= e \vee f$. This nine-element partial algebra has $\sigma_5(U_1)=21.43$. This case excluded because $\sigma_5(U_1)=21.43<25$.
	
    The next case is  $t=8$. Then  by the symmetry $z$ is not new element, either $a$, $c$, $x$ or $y$ recalling   that $a$ and $b$ are symmetric and $c$ and $d$ are symmetric. We will consider the following cases.
	 
	\bf Case i: \rm If $z= e \vee f=a$, this case is covered by taking the partial algebra  $U_2=\{a,b,c,d,e,f,x,y\}$ with  $x= a \vee b$, $y= c \vee d$ and $a= e \vee f$. It has  $\sigma_5(U_2)=21$ using the computer program. Like the above, it is ruled out because  $\sigma_5(U_2)=21 <25$. 
	
	\bf Case ii: \rm If $z= e \vee f=c$, this case captured by taking the partial algebra  $U_3=\{a,b,c,d,e,f,x,y\}$ with  $x= a \vee b$, $y= c \vee d$ and $c= e \vee f$. It has  $\sigma_5(U_3)=21$ using the computer program. Like above it is excluded because  $\sigma_5(U_3)=21 <25$. 
	
	\bf Case iii: \rm  If $z= e \vee f=y$, this case covered by taking the partial algebra  $U_4=\{a,b,c,d,e,f,x,y\}$ with  $x= a \vee b$, $y= c \vee d$ and $y= e \vee f$. It has  $\sigma_5(U_4)=21.8$ using the computer program. Like above it is ruled out because  $\sigma_5(U_4)=21.8 <25$. 
	
	\bf Case iv: \rm  If $z= e \vee f=x$, this case captured by taking the partial algebra  $U_5=\{a,b,c,d,e,f,x,y\}$ with  $x= a \vee b$, $y= c \vee d$ and $x= e \vee f$. It has  $\sigma_5(U_5)=21.8$ using the computer program. Like above it is excluded because  $\sigma_5(U_5)=21.8 <25$. 
	
	The case where $t=7$. This means neither $y$ nor $z$ is a new element, recall that as $a$ and $b$ are symmetric, one can assume that either  $y= c \vee d=a$, $z= e \vee f=b$ or $y= c \vee d= z = e \vee f=x$. now let us consider the first case. This case is covered by taking the partial algebra  $U_6=\{a,b,c,d,e,f,x\}$ with  $x= a \vee b$, $y= c \vee d=a$ and $z= e \vee f=b$.  It has  $\sigma_5(U_6)=20.5$ using the computer program. Like the above, it is ruled out because  $\sigma_5(U_6)=20.5 <25$. For the case that $y= c \vee d= z = e \vee f=x$, it is covered by taking the partial algebra  $U_7=\{a,b,c,d,e,f,x\}$. It has  $\sigma_5(U_7)=22.75 <25$.
	
	For the case where $t=6$,   $x= a \vee b$ is not a new element and cannot be any of the $a,b,c,d,e,f$ because we assumed that all the elements were incomparable, which contradicts $t=6$. So this case is ruled out.

    Note that there is a non-comparable  $\{b,d\}$ disjoint from $\{c,d\}$ with $b<c$, and we have another element $a$, there are cases that depend on the position of $a$,  so to complete the proof, consider the following cases.
    
	\bf Case i: \rm Here, $a$ is incomparable with $d$ and $b$ and so $c$ because $b<c$. By using the computer program the join semilattice here
     has the relative number of subuniverses equal to 22.
 
	
	\bf Case ii: \rm Here, $a$ is comparable with $d$ and incomparable with $b$ and hence with $c$. this case is covered by taking the partial algebra  $H=\{a,b,c,d,1\}$ with  $ 1= a \vee b$, $1= c \vee a$. It has  $\sigma_5(H)=21$ using the computer program.  It is ruled out because  $\sigma_5(H)=23 <25$.

		
	\bf Case iii: \rm Here, $a$ is comparable with $b$ and then $a<b<c$, and it is  incomparable with $d$. The join-semilattice is $U_{10}=\set{a,b,c,d,1}$ with $ c\vee d= a\vee d=b\vee d=1$. It has $\sigma_5(U_{10}) =25,$ this subsemilattice (isomorphic to) $H_5$. 
	Now that all the other possibilities have been excluded, we know that $ H_5 $ is a join-subsemilattice of $(L,\vee)$. Observe that $ H_5 $ has no narrows. Therefore, by Lemma \ref{lemma2.1-7}\ref{lemm21c}, $(L,\vee)$ is of the desired form.

 Now, we will prove part \ref{thmmaine}.
 
	Assume that $ (L,\vee)$ is of the given form. Then $\sigma_5(L,\vee)=24.5$ is clear from Lemma \ref{lemma2.1-7}\ref{lemm21c} and Lemma \ref{lemmaLHFtfsZh}\ref{iii}. To prove the converse, assume that $\sigma_5(L,\vee)=24.5$. For this, we will not consider the cases that are excluded in Theorem 1 in \cite{ahmed2021first} or in the proof above. We will divide the proof in two cases depending on the number of anti-chains. 
 
 \bf Case i: \rm If we have two 2-element antichain $\{a,b\}$ and $\{c,d\}$. Now let $x:=a\vee b$ and $y:=c\vee d$. There are sub-cases depending on the position of $y$. Some of this sub-cases will excluded by the symmetric of $a$ and $b$. 
 
 \bf Sub-case i: \rm This case is covered by taking the partial algebra $U_{11}=\{a,b,c,d,x,y\}$ with edges $cy,dy,yb,bx,ax$. It has $\sigma_5(U_{11})=22.5$ by using the computer program. 
 
 \bf Sub-case ii: \rm This case is covered by taking the partial algebra $U_{12}=\{a,b,c,d,x,y\}$ with edges $cy,dy,yb,ya,bx,ax$. A part from the notation used, this join subsemilattice is isomorphic to $H_3\gluplus B_4$. It has $\sigma_5(U_{12})=24.5$ by Lemma \ref{lemmaLHFtfsZh}~\ref{ii}.
 
 \bf Case ii: \rm If we have three 2-element antichains $\{a,b\}$, $\{c,d\}$ and $\{e,f\}$ with $x:=a\vee b$, $y:=c\vee d$ and $z:=e\vee f$. Depending on the position of $y$ and $z$ we will have cases as follows. The number of possible sub-cases can be reduced by symmetric role $a$ and $b$ or $c$ and $d$ play. 

  \bf Sub-case i: \rm The partial algebra $U_{13}=\{a,b,c,d,e,f,x,y,z\}$ with edges\newline $cy, dy, yb, bx, ax, ez, fz$, has $\sigma_5(U_{13})=19.6875$ by using the computer program. 
  
  \bf Sub-case ii: \rm The partial algebra $U_{14}=\{a,b,c,d,e,f,x,y,z\}$ with edges \newline $cy, dy, yb, ya, bx, ax, ez, fz$, has $\sigma_5(U_{14})=21.4375$ by using the computer program. 
  
  \bf Sub-case iii: \rm The partial algebra $U_{15}=\{a,b,c,d,e,f,x,y,z\}$ with edges \newline   $cy, dy, yb, za, bx, ax, ez, fz$, has $\sigma_5(U_{15})=16.9375$ by using the computer program. 
  
  \bf Sub-case iv: \rm Let assume that $z$ is not a new element, let $z=c= e\vee f$. This case is covered by taking the partial algebra $U_{16}=\{a,b,c,d,e,f,x,y\}$ with edges $cy, dy, yb, bx, ax, ec, fc$. It has $\sigma_5(U_{16})=18.75$ using the computer program. 
  
  \bf Sub-case v: \rm The partial algebra $U_{17}=\{a,b,c,d,e,f,x,y,z\}$ with edges \newline $cy, dy, yb, zc, bx, ax, ez, fz$, has $\sigma_5(U_{17})=18.94$ by using the computer program. 
  
  \bf Sub-case vi: \rm This case is covered by taking the partial algebra $U_{18}=\{a,b,c,d,e,f,x,y,z\}$ with edges $cy, dy, yb, zc, zd, bx, ax, ez, fz$. It has $\sigma_5(U_{18})=18.9375$ by using the computer program. 
  
  \bf Sub-case vii: \rm The partial algebra $U_{19}=\{a,b,c,d,e,f,x,y,z\}$ with edges \newline $cy, dy, yb, ya, zc, zd, bx, ax, ez, fz$, has $\sigma_5(U_{19})=21.4375$ by using the computer program.

Finally, we prove part \ref{thmmainf}.

let $(L,\vee)$ be an n-element semilattice. We know from Lemma\ref{lemma2.1-7}(iii) and Lemma \ref{lemmaLHFtfsZh}\ref{iiii} that if 
\begin{center}
 $(L,\vee)\cong K_3 \gluplus C_1$ or
		$(L,\vee)\cong C_0 +_{ord} K_3 \gluplus C_1$, 
  \end{center}
where \badgood{$C_0$ , $C_1$ and  $C_2$ }{$C_0$ and $C_1$ } are finite chains, then $\sigma_5(L,\vee) = \sigma_5(K_3)=24$ indeed. Conversely, assume that $\sigma_5(L,\vee) = 24$ and it has three incomparable elements $a,b,c$. Clearly, $\{a,b,c,a\vee b\vee c\}$ is a join-semilattice isomorphic to $K_3$. But $\sigma_5(K_3)=24$ by Lemma\ref{lemma2.1-7}(iii) and Lemma \ref{lemmaLHFtfsZh}\ref{iiii} immediately tell us that $(L,\vee)$ is of the desired form. 

 From this, the proof of Theorem \ref{thmmain} is complete.

\end{proof}


\bibliographystyle{unsrt}
\bibliography{ref}

\end{document}